\documentclass[12pt]{amsart}

\theoremstyle{definition}

\theoremstyle{parrafo}

\begin{document}

\title[]{The  weak type $(1,1)$ bounds for the maximal function associated to cubes
grow to infinity with the dimension}

\author{J. M. Aldaz}
\address{Departamento de Matem\'aticas,
Universidad  Aut\'onoma de Madrid, Cantoblanco 28049, Madrid, Spain.}
\email{jesus.munarriz@uam.es}


\thanks{The author was partially supported by Grant MTM2009-12740-C03-03
 of the
D.G.I. of Spain}







\begin{abstract} Let $M_d$ be the centered Hardy-Littlewood maximal function
associated to cubes in $\mathbb{R}^d$  with Lebesgue measure, and
let  $c_d$ denote the lowest constant 
appearing
in the weak type (1,1) inequality
satisfied by $M_d$.
 We show that 
$c_d \to \infty$
as $d\to \infty$, thus answering, for the case of cubes, a long
standing open question of E. M. Stein
and J. O. Str\"{o}mberg.
\end{abstract}


\maketitle


\section {Introduction and result.}

\markboth{J. M. Aldaz}{Weak type $(1,1)$ bounds}

By a cube $Q(x,r)$ we mean a closed $\ell_\infty$ ball of radius $r$
and center $x$ in $\mathbb{R}^d$, that is,
a closed cube centered at $x$, with sides parallel to the coordinate axes, and 
sidelength $2r$.
Let $M_d$ be the centered maximal function 
\begin{equation}\label{HLMF}
M_{d} f(x) := \sup _{r > 0} \frac{1}{|Q(x, r)|} \int _{Q(x, r)} \vert f(y)\vert dy
\end{equation}
associated to cubes and  Lebesgue measure in $\mathbb{R}^d$.
A fundamental feature of the Hardy-Littlewood maximal function $M$
is that it satisfies the weak-type $(1,1)$ inequality: There exists
a constant $c>0$ such that for all $\alpha > 0$ and all $f\in L^1$,
\begin{equation}\label{weaktype} 
\alpha |\{ Mf \ge \alpha\}| \le c \|f\|_1.
\end{equation}
Denote by $c_d$ the best (i.e. lowest) constant 
satisfying (\ref{weaktype}) in $\mathbb{R}^d$. 
\vskip .3 cm

\noindent{\bf Theorem.} {\em  
Fix $T > 0$. Then there exists a $D = D(T)$ such that for every
dimension $d \ge D$, $c_d \ge T$.}

\vskip .3 cm

Thus, $c_d \to \infty$
as $d\to \infty$. In fact, these constants approach $\infty$ in a
monotone manner, since $c_{d+1} \ge c_d$ by \cite[Theorem 2]{AV}.

It is well known that given $1 < p \le \infty$, there exists
a constant $c_p$ such that for all $f\in L^p (\mathbb R^d)$, $\|Mf\|_p\le c_p \|f\|_p$.
When $p =\infty$, trivially $c_p = 1$ in every dimension, for averages
never exceed a supremum. Dimension independent estimates are useful
whenever one is interested in extending results from the finite dimensional to the infinite dimensional
setting.
For the
maximal function associated to euclidean balls,  E. M. Stein showed
 that one can take $c_p$ to be independent of $d$ (\cite{St1}, \cite{St2}, \cite{StSt}, see also \cite{St3}).
Stein's result was
generalized to the maximal function defined using balls given by  arbitrary norms by J. Bourgain (\cite{Bou1},
\cite{Bou2}, \cite{Bou3}) and A. Carbery (\cite{C}) when $p>3/2$.
Given $1\le q <\infty$, the $\ell_q$ balls
are defined using the norm $\|x\|_q :=\left( |x_1|^q+ |x_2|^q+\dots + |x_d|^q\right)^{1/q}$, and the $\ell_\infty$ balls,
 using  $\|x\|_\infty :=\max_{1 \le i \le d} \left\{ |x_1|,  |x_2|, \dots , |x_d|\right\}$.
For $\ell_q$ balls, $1\le q <\infty$, D. M\"{u}ller \cite{Mu}
showed that uniform bounds again hold for every $p
> 1$. 
With respect to weak type bounds, in \cite{StSt} E. M. Stein
and J. O. Str\"{o}mberg proved  that the  
best
constants in the weak type (1,1) inequality satisfied by  the maximal
function associated to arbitrary balls  grow
at most like $O( d \log d)$, while if the balls are euclidean, then 
the best constants grow at most like
$O(d)$. They also asked if uniform bounds could be found. The
theorem above shows that in the case of cubes the
answer is negative. If the $d$-dimensional Lebesgue measures are replaced by a sequence of finite,
absolutely continuous radial measures with decreasing densities (such as,
for instance,
the standard Gaussian
measures) then best constants grow exponentially with $d$, cf. \cite{A2}.

In recent years evidence has been mounting
to the effect that 
not only weak type $(1,1)$ inequalities are formally stronger
than strong $(p,p)$ inequalities for $1 < p < \infty$  (since
the latter are implied by the former via interpolation) but
they are also stronger in a substantial way, meaning that
the strong type may hold for all $p > 1$ while the weak type $(1,1)$
may fail. This is the case, for instance, with the uncentered maximal
function associated to the standard gaussian measure
and euclidean balls. It is shown in
\cite{Sj} that this maximal function is not of weak type $(1,1)$,
while it is strong $(p,p)$ for all $p >1$, cf. \cite{FSSjU} (for cubes
the strong $(p,p)$ type follows from a more general result in
\cite[cf. Theorem 1]{CF}). The theorem above may represent another instance of this phenomenon,
with respect to uniform bounds in $d$. However, it is not
known for cubes whether uniform bounds hold when $1 < p \le 3/2$
(it is suggested in \cite{Mu} that the answer may be negative, and conjectured
in \cite{ACP} that the answer is positive).

Before presenting the proof, we make some comments on the method 
of
 discretization for weak type $(1,1)$ inequalities. It consists in replacing $L^1$ functions by
 finite sums of Dirac deltas. 
 This leads to elementary arguments of a combinatorial nature.
  The fact that one can get lower bounds for
 $c_d$ 
 using Dirac deltas instead of functions is obvious, by mollification. And this all we need here.
 
We mention for completeness that considering Dirac deltas also suffices to give upper bounds, as shown  by M. de Guzm\'an, 
see \cite[Theorem 4.1.1]{Gu}. Furthermore,  
M. Trinidad Men\'arguez and F. Soria proved that discretizing does not
alter constants, cf.  \cite[Theorem 1]{MS}, so it can be used to
study the precise values of $c_d$. This method was utilized,   for instance, in  \cite{A},
were it was shown that $37/ 24 \le c_1\le  \frac {9 + \sqrt{41}}8$, thereby 
refuting the conjecture that  $c_1=3/2$ (cf. \cite[Problem 7.74 c]{BH}) and 
showing that $c_1 < 2$, which is the best constant in the uncentered case. 
Discretization was 
also used in
\cite{Me1} and  \cite{Me2}, where 
the exact value of $c_1=\frac {11 + \sqrt{61}}{12}$ was found. No
best constants are known for dimensions larger than one.

Let us point out that the configuration of Dirac deltas we will 
utilize
had previously been considered in 
\cite[Theorem 6]{MS}, for the same purpose of bounding
$c_d$ from below. It is shown there that 
$c_d \ge \left( \frac{1+ 2^{1/d}}{2} \right)^d$. 
But this
yields no information as to whether there is a uniform 
upper bound for $c_d$, since $\left( \frac{1+ 2^{1/d}}{2} \right)^d < 2$. 

The first version of this article contained a simple counting error at the
beginning of the proof, which  rendered large parts of it useless. I am most indebted to professor Keith Rogers for 
pointing out this mistake to me;  substantial modifications of the
argument were required in order to fix it. 
I am also indebted to
professors Javier P\'erez L\'azaro and Peter Sj\"ogren, and to an anonymous
referee, for carefully reading this paper and for making several useful suggestions, which led to a thorough rewriting and simplification of 
this paper.

\section
{Proof.}

Given a locally finite measure $\nu$ in
$\mathbb{R}^d$, the maximal function $M_d\nu$ 
is defined by 
\begin{equation}\label{maxmeas}
M_d\nu (x) := \sup _{r > 0}  \frac{\nu Q(x, r)}{|Q(x, r)|}.
\end{equation}

For notational simplicity we start considering the infinite measure $\mu^d$ in
$\mathbb{R}^d$ obtained by placing one Dirac delta at each point
of the integer lattice $\mathbb{Z}^d$. The finite measure
exhibiting a lower bound for $c_d$ will then be obtained by restricting $\mu^d$
to a sufficiently large cube. Note that 
$\mu^d = \mu^1 \times \mu^1 \times\cdots\times \mu^1$. At first, we
will work within the unit cube $[0,1]^d$ only.

Given $u\in (0, 1)$ and an interval $I\subset \mathbb{R}$, call $y\in I$
centered at level $u$ (more briefly, centered, or  $u$ centered)
if it belongs to the closed subinterval with the same center
and   length $(1 - u) |I|$, and off center (at level $u$) otherwise. In particular,
for $I = [0,1]$ the centered points are those in $[2^{-1} u, 1 - 2 ^{-1}u]$.
The role of $u$ in the proof is to serve as a discrete parameter,
used to describe which cubes should be considered when estimating the
value of the maximal function at a given point.

It can be shown that the maximal function is large on the set 
$E_u\subset [0,1]^d$
of points $x = (x_1, \dots,x_d)$ with many centered coordinates, where
``large" is determined by a fixed $t >> 1$, and 
``many"  means more than $(1 - u) d + t\sqrt{d u(1-u)}$. Since 
$t\sqrt{d u(1-u)}$ 
amounts to  $t$ standard deviations 
of a binomially distributed random variable with parameters
$d$ and $u$, the Central Limit Theorem allows us to bound $|E_u|$ 
from below, provided $d$ is large enough (we mention that a similar
argument can be carried out on the set 
of points with many uncentered coordinates). For a fixed $u$, the measure of $E_u$ 
 as $d\to \infty$ turns out to be too small, since we are $t >> 1$ standard
 deviations away from the mean.
  On the other hand, estimates
for the size of $M_d\mu^d$ worsen when we have roughly $(1 - u) d$
 centered  coordinates. Changing the value of $u$ by discrete
 steps and taking the union of many  $E_u$'s, we obtain a sufficiently large set over which $M_d\mu$
can be shown to take high values (unlike $u$, the value of $t$ is fixed throughout
the argument, so dependency on $t$ is not indicated in the notation).
In order to control the intersections of different $E_u$'s, it is useful to
also bound from below the number of uncentered coordinates, so actually
we shall slightly modify these sets and call them $E^u$ instead
of $E_u$.

\vskip .2 cm

Fix $t >> 1$. The assumption that $t$ is very large will be used without
further mention (save for some occasional reminder). 
But we emphasize that the value of $t$ remains unchanged throughout
the proof; in particular, it does not approach $\infty$ as $d\to\infty$.
So we  will assume, again without
further mention, that expressions such as $t/\sqrt{d}$ are as small as
needed  each time they appear. 
 
Recall that the standard deviation of a Bernoulli  trial with parameter
$u$ is $\sigma_u = \sqrt{u(1-u)}$.  Define, for each $u\in [1/8,1/4]$, 
\begin{equation}\label{setE}
E^{u}:=\{x\in [0,1]^d: \mbox{ the number } k \mbox{ of coordinates }
j_1,\dots, j_k \mbox{ for which } 
\end{equation}
\begin{equation*}
x_{j_i}\in \left[0, 2^{-1}u\right)\cup 
\left(1 - 2^{-1} u, 1\right] \mbox{ satisfies }
u d -   (t + t^{-1}) \sigma_u \sqrt{d} < k \le u d - 
t \sigma_u \sqrt{d}\}.
\end{equation*}
 The values $1/8$ and $1/4$ are of no special significance; we could have fixed any $0 <a < b < 1$ and chosen
$u\in [a,b]$ instead.

In order to prove the theorem,  first we estimate the size of $E^u$ for each
$u = 1/8, 1/8 + t^{- 4/3}, 1/8 +  2 t^{- 4/3}, 1/8 +  3 t^{- 4/3}\dots \le 1/4$, so 
we consider $\Theta( t^{4/3})$
 different values of $u$, where, as usual, $\Theta$ stands for exact order.
Second, using the fact that distinct values
of $u$ differ by at least $t^{- 4/3}$,  
we prove that
different sets $E^u$ have very small intersection.
Third, we take
the union of the $\Theta( t^{4/3})$ sets $E^u$ and bound the measure
of this union from below. 
And fourth, we show that $M\mu^d$ is large on
each $E^u$, and hence on their union.

Up to here, the argument is carried inside
$[0,1]^d$. To complete the proof, we replace $\mu^d$ by a finite measure,
and apply the estimates obtained within  $[0,1]^d$ to several translates of it.

Let $Z\sim N(0,1)$ denote a standard normally distributed random variable.

\vskip .2 cm

{\em Claim 1: For all $u\in [1/8, 1/4]$ there exists
a $D = D(u)$ such that if
$d \ge D$, then  
\begin{equation} \label{eubounds}
\frac{e^{-t^2/2} }{2 e^2 t\sqrt{2 \pi}}  
< 
|E^{u}| <
\frac{e^{-t^2/2}}{t\sqrt{\pi}}. 
\end{equation} 
}

\vskip .2 cm

{\em Proof.} Define a collection of independent Bernoulli random variables $X_{u,i}$
by setting $X_{u,i}(x) = 1$ 
if the $i$-th coordinate of 
 $x\in[0,1]^d$ satisfies $x_i \in \left[0, 2^{-1}u\right)\cup 
\left(1 - 2^{-1} u, 1\right]$, and $X_{u,i}(x) = 0$ otherwise.
Then the probability of having exactly
 $k$ off center and $d - k$ centered coordinates is 
\begin{equation*}
\binom{d}{k} u^{k} (1 - u)^{d-k}.
\end{equation*}
Set $S_{u,d}:= \sum_{i = 1}^d X_{u,i}$, so $S_{u,d}$
counts the number of uncentered coordinates. Then $S_{u,d}\sim B( u, d)$ is  binomially distributed
with  mean $E(S_{u,d})= u d$ and standard deviation
$\sigma (S_{u,d}) = \sqrt{d u(1 - u)} = \sigma_u \sqrt{d}$,
where $\sigma_u$ is the standard deviation of $X_{u,i}$.
Thus, the Lebesgue measure of $E^u$ is given
by 
$$
|E^u| =  P \left(u d - 
(t + t^{-1}) \sigma_u  \sqrt{d} < S_{u,d} \le u d - 
t\sigma_u \sqrt{d}\right).
$$ 
 By the Central Limit Theorem,  for all 
$d$ large enough we have
\begin{equation}\label{sdn}
2^{-1} P(-t - t^{-1} < Z \le - t) < |E^u| < 
\sqrt{2} P(-t - t^{-1} < Z \le - t).
\end{equation}
Since 
$ P(-t - t^{-1} < Z \le - t ) = 
\frac{1}{\sqrt{2 \pi}} \int_{t}^{t + t^{-1}} e^{-y^2/2} dy$
and $e^{-2} e^{-t^2/2} \le e^{-y^2/2} \le e^{-t^2/2}$ for every
$y \in [t, t + t^{-1}]$, we obtain
(\ref{eubounds}).
\qed

\vskip .2 cm

We show next that if $u$ and $v$ are ``far apart", then 
 $|E^{v}\cap E^{u}|$ is small relative to $|E^u|$. 
 
\vskip .2 cm

{\em Claim 2: Fix  $u, v\in [1/8, 1/4]$ with $u - v\ge t^{-4/3}$. 
Then there exists a $D = D(u)$ such that for all $d\ge D$,  
we have 
\begin{equation} \label{euvbounds}
|E^{v}\cap E^{u}|\le t^{-1/3} e^{-2 t^{2/3}/9} |E^u|.
\end{equation}}

\vskip .2 cm

{\em Proof.} We partition the sets $E^u$ into subsets $A_{u, K}$,
consisting of all points with coordinates $x_j$ off center 
if and only if
$j\in K$. More precisely, 
let us fix a subset
  $K\subset\{1, \dots, d\}$ of cardinality $k$, with $k$ 
satisfying  
\begin{equation} \label{fixk} 
u d - 
(t + t^{-1}) \sigma_u \sqrt{d} < 
k \le u d - 
t\sigma_u \sqrt{d}.
\end{equation}
We define
\begin{equation}\label{setA}
A_{u, K}:=\left\{x\in [0,1]^d: x_j \in \left[0, 2^{-1}u\right)\cup 
\left(1 - 2^{-1} u, 1\right] \mbox{ if and only if } j \in K\right\}.
\end{equation}
To check that $A_{u, K} \cap A_{u, M} = \emptyset$ unless $K = M$,
suppose  that there is a $j\in (K\setminus M) \cup (M\setminus K)$.
Then for any $x\in A_{u, K} \cap A_{u, M}$, its $j$-th coordinate
must be simultaneously centered and off center at level $u$, which
is impossible. Since $E^u = \cup_K A_{u, K}$, we do have a partition of $E^u$.

Let $v \le u - t^{-4/3}$. In order to estimate $|A_{u,K}\cap E^v|$, consider an arbitrary  set
$A_{v,M}$ in the partition of $E^v$. We may suppose that $M\subset K$,
since otherwise $A_{u,K} \cap A_{v,M} = \emptyset$, for the same reason
as before: If there is a $j\in M\setminus K$ and 
$x\in A_{u,K} \cap A_{v,M}$, then
 $x_j$
must be $u$ centered and $v$ off center.
Suppose $M$ has cardinality $m$, and let $x \in A_{u,K} \cap A_{v,M}$. Observe that 
$x_i\in [0, 2^{-1} v)\cup (1- 2^{-1} v, 1]$
 for every
 $i\in M$, $x_j\in [2^{-1} v,  2^{-1} u) \cup (1- 2^{-1} u, 1 - 2^{-1} v]$ for 
 every
$j\in K\setminus M$, and $x_r\in [2^{-1} u, 1- 2^{-1} u]$
 for the
remaining $d - k$ coordinates.
Thus
$$
|A_{u,K} \cap A_{v,M}| = v^m \times (u - v)^{k -m}\times (1 - u)^{d -k}
= u^k (1 - u)^{d -k} \left(\frac{v}{u}\right)^m \left(1 -\frac{v}{u}\right)^{k - m}.
$$
The lower
bound on $k$ given by (\ref{fixk}) allows us to conclude that 
 $k > m$ for sufficiently large $d$, since in that case
 $u d - 
(t + t^{-1})\sigma_u \sqrt{d} > v d - 
t\sigma_v \sqrt{d}$. Only the upper bound $m\le [v d - t \sigma_v \sqrt{d}]$ (where $[w]$ denotes the integer part of $w$) is needed in the next estimate. 
 Summing first over 
all sets  $M\subset K$ of fixed cardinality $m$, and then over all $m\le [v d - t \sigma_v \sqrt{d}]$, we get
\begin{equation}\label{aka}
|A_{u,K} \cap E^{v}| \le u^k (1 - u)^{d -k} 
\sum_{m= 0}^{[v d - 
t\sigma_v \sqrt{d}]}
\binom{k}{m}\left(\frac{v}{u}\right)^m \left(1 -\frac{v}{u}\right)^{k - m}.
\end{equation}
As before, we control the sum above by using the
Central Limit Theorem, applied to the binomially distributed random 
variable $S_{v/u, k} \sim B( v/u, k)$. We shall need a lower bound for 
$E (S_{v/u,k})$ and an upper bound for $\sigma(S_{v/u,k})$.
From $u d - 
(t + t^{-1}) \sigma_u \sqrt{d} < 
k \le u d$ we obtain  
\begin{equation}\label{Ek}
E (S_{v/u,k})  = k\left(\frac{v}{u}\right)
>
v d - 
\left(t + t^{-1}\right) v \sqrt{d}\sqrt{u^{-1} -1}
\end{equation}
and
\begin{equation}\label{sigmak}
\sigma(S_{v/u,k}) = \sqrt{k \left(\frac{v}{u}\right)\left(1 - \frac{v}{u}\right)} < \sqrt{ v d   \left(1 - \frac{v}{u}\right)}.
\end{equation}
Since $t >> 1$,  $2 t^{-2/3} \le \sqrt{1 - v/u}$, and $\sqrt{v} \sqrt{u^{-1}-1} + 
\sqrt{1 -v}<\sqrt{7}$, 
\begin{equation*}
\frac{2 t^{1/3}}{3} 
<
\frac{t \sqrt{1 - \frac{v}{u}}}{\sqrt{7}}
 - \frac{\sqrt{7}}{t^{1/3}}
<
\frac{t \left(1 - \frac{v}{u}\right)}{\left(\sqrt{v} \sqrt{\frac{1}{u}-1} + 
\sqrt{1 -v}\right)\sqrt{1 - \frac{v}{u}}}
  - \frac{\sqrt{v} \sqrt{\frac{1}{u}-1}}{t \sqrt{1 - \frac{v}{u}}}
 \end{equation*}
\begin{equation*}
=
\frac{t \left(\sqrt{1 -v} - \sqrt{v} \sqrt{\frac{1}{u}-1}\right) - t^{-1}  \sqrt{v} \sqrt{\frac{1}{u}-1}}{\sqrt{1 - \frac{v}{u}}}
\end{equation*}
\begin{equation*}
=
\frac{t\sqrt{vd}\sqrt{1 -v} - 
\left(t + t^{-1}\right) v \sqrt{d} \sqrt{\frac{1}{u}-1} + vd -vd}{\sqrt{vd}\sqrt{1 - \frac{v}{u}}}
<
\frac{E (S_{v/u,k}) - v d + 
t\sigma_v \sqrt{d}}{\sigma(S_{v/u,k})},
\end{equation*}
where the last inequality follows from (\ref{Ek}) and (\ref{sigmak}).
Hence, by the Central Limit Theorem we have, for all sufficiently large $d$, 
\begin{equation}\label{skn}
P \left(S_{v/u,k}  \le v d - 
t\sigma_v \sqrt{d} \right) \le P \left(\frac{S_{v/u,k} - 
E (S_{v/u,k})}{\sigma(S_{v/u,k})}  \le - \frac{2 t^{1/3}}{3} \right)
\end{equation}
\begin{equation}\label{skn2}
\le 
\frac{2 \sqrt{2 \pi}}{3} P\left(Z \le - \frac{2 t^{1/3}}{3}\right)
=  
\frac{2}{3} \int_{\frac{2 t^{1/3}}{3}}^{\infty} e^{-y^2/2} dy
\le  
\frac{1}{t^{1/3}} \int_{\frac{2t^{1/3}}{3}}^{\infty} y e^{-y^2/2} dy 
= 
\frac{e^{-2 t^{2/3}/9}}{t^{1/3}}.
\end{equation}
Thus, from (\ref{aka}) and (\ref{skn})-(\ref{skn2}) we get 
\begin{equation}\label{akad}
|A_{u,K} \cap E^{v}| \le u^k (1 - u)^{d -k} 
t^{-1/3}  e^{-2 t^{2/3}/9},
\end{equation}
and now (\ref{euvbounds}) follows by adding up over all the sets
$A_{u, K}$ in the partition of $E^u$.
\qed

\vskip .2 cm

Next, we  write $u(j) := 1/8 + j t^{- 4/3}$, letting $u$ range over
$[1/8, 1/4]$ in discrete steps of size $t^{- 4/3}$.

\vskip .2 cm

{\em Claim 3: Let $j= 0, 1, 2, \dots M$, where $M$ is the largest 
integer $j$   satisfying $j t^{- 4/3} \le 1/8$. Then 
\begin{equation}\label{levelset}
\left|\bigcup_{j=0}^M E^{u(j)}\right| 
 > \frac{t^{1/3} e^{-t^2/2}}{20 e^2 \sqrt{2 \pi}}. 
 \end{equation}}

\vskip .2 cm

{\em Proof.} Let $0 \le k < j\le M$ be natural numbers. We apply Claim 2 
and Claim 1 to all pairs  $u(k) < u (j)$, obtaining
\begin{equation}\label{intersec}
\sum_{0\le k < j \le M} \left| E^{u(k)} \cap E^{u(j)}\right| < t^{8/3} \left(\frac{e^{-2 t^{2/3}/9}}{t^{1/3}} \right)
\left( \frac{e^{-t^2/2}}{t}\right) 
= t^{4/3} e^{- 2 t^{2/3}/9} e^{-t^2/2} = O\left(t^{-1} e^{-t^2/2}\right).
\end{equation}
Using the inclusion exclusion principle, together with (\ref{intersec}) and the lower bound from (\ref{eubounds}),
we obtain 
\begin{equation}\label{unionlb}
\left|\bigcup_{j=0}^M E^{u(j)}\right| 
\ge
\sum_{j=0}^M \left|E^{u(j)}\right| 
-
\sum_{0\le k < j \le M} \left| E^{u(k)} \cap E^{u(j)}\right|
\end{equation}
\begin{equation} 
>
\left(\frac{t^{4/3}}{8}\right) \frac{e^{-t^2/2} }{2 e^2 t\sqrt{2 \pi}} - 
O\left(t^{-1} e^{-t^2/2}\right) > \frac{t^{1/3} e^{-t^2/2}}{20 e^2 \sqrt{2 \pi}}. 
\end{equation}
\qed

\vskip .2 cm

After having estimated the measure of
$\bigcup_{j=1}^M E^{u(j)}$, it remains to show that $M\mu^d$
takes high values on this set. We do this
next.
 
\vskip .2 cm

{\em Claim 4: Fix  $u\in [1/8, 1/4]$. Then 
$E^{u}\subset \{M\mu^d > 2^{-1} e^{t^2/2}\}$.}

\vskip .2 cm

{\em Proof.} Let $x_j\in [2 ^{-1} u, 1- 2 ^{-1} u]$. Given any integer $s > 0$, 
$\mu^1 [x_j - (s -2 ^{-1} u), x_j + s - 2 ^{-1}  u] = 2 s$. Suppose 
$y\in [0,1]$ is off center, say for instance $y > 1 - 2 ^{-1}  u$. Shifting
the  interval $[x_j - (s -2 ^{-1} u), x_j + s - 2 ^{-1}  u]$ to the right by $y - x_j$ (so now it is centered at $y$) loses 
at most one Dirac delta on the left. Thus,  
$\mu^1([y - (s - 2 ^{-1}  u), y + s - 2 ^{-1}  u]) \ge 2 s - 1$. 
Suppose $x\in [0,1]^d$ has $r$ off center and $d-r$ centered coordinates,
where  $r \le r_0 := u d - t\sigma_u \sqrt{d}$.
Then for every $s=1, 2, 3, \dots$,
\begin{equation}\label{maxlevel}
M_d\mu^d (x) \ge 
\frac{(2s)^{d-r}(2s - 1)^r}{\left(2\left(s - 2 ^{-1} u \right)\right)^d} = 
\frac{\left(1 - \frac{1}{2 s}\right)^r}{\left(1 - \frac{u}{2 s}\right)^d}
\ge \frac{\left(1 - \frac{1}{2 s}\right)^{r_0}}
{\left(1 - \frac{u}{2s}\right)^d}.
\end{equation}
The next step consists in showing that for  some
suitably chosen $s$, the right hand side of (\ref{maxlevel}) is bounded
below by $2 ^{-1} e^{t^2/2}$.
Set 
$f(s):= \left(1 - \frac{1}{2 s}\right)^{r_0}
/\left(1 - \frac{u}{2 s}\right)^d$.
An elementary calculus argument shows that  $f(s)$ 
is maximized over $s \ge 1$ when 
\begin{equation}\label{maxs}
s =\frac{u \left(d - r_0\right)}{2\left(u d - r_0 \right)}=: s_0,
\end{equation}
and this is the only critical point, so $f$ decreases as we move away
from $s_0$. In particular, $f(s_0) \ge f([s_0]) \ge f(s_0 -1)$
(where $[s_0]$ denotes the integer part of $s_0$) 
so for convenience we shall use $s_0 -1$ instead of $[s_0]$ in (\ref{maxlevel}).
 Thus, 
\begin{equation}\label{maxfnlevel}
\log M_d\mu^d (x) \ge 
 \left(u d - t\sigma_u \sqrt{d}\right) \log \left(1 -  
  \frac{1}{2(s_0 - 1)}\right) - d \log \left(1 -  
  \frac{u}{2(s_0 - 1)}\right).
\end{equation}
Replacing $r_0$ by its value in (\ref{maxs}) we see that  
$$
s_0 =  \frac{\sigma_u \sqrt{d}}{2 t} + \frac{u}{2} = \Theta (\sqrt{d})
\mbox{ \ \ \ and \ \ \ } 2(s_0 - 1) = \frac{\sigma_u \sqrt{d}}{t}\left(1 - \frac{(2 - u) t}{\sigma_u\sqrt{d}}\right).
$$
Thus, using $(1-w)^{-1} = 1 + w + O(w^2)$ applied to $w_1 := \frac{(2 - u) t}{\sigma_u\sqrt{d}} = \Theta (1/\sqrt{d})$ we get
\begin{equation}\label{w_2}
  \frac{1}{2(s_0 - 1)} = \frac{t}{\sigma_u \sqrt{d}} + \frac{(2 - u) t^2}{\sigma_u^2 d} + O\left(\frac{1}{d^{3/2}}\right).
\end{equation}
Finally, from $\log(1-w) = - w - w^2/2 + O(w^3)$ applied to 
$w_2 := \frac{1}{2(s_0 - 1)} = \Theta (1/\sqrt{d})$
and  to $w_3 := u w_2$, we obtain, by substituting in (\ref{maxfnlevel})
and simplifying,
\begin{equation}\label{maxfnlevel2}
\log M_d\mu^d (x) \ge   \frac{t^2}{2} + O\left(  
  \frac{1}{\sqrt{d}}\right).
\end{equation}
 \qed

\vskip .2 cm

{\em   Completing the argument.} The last step consists in fixing
$d$ (so large that all the preceding estimates hold) and replacing the
infinite measure $\mu^d$ by a finite measure $\mu^d_R$, such that the ratio of
unit volume cubes to Dirac deltas is close to  1. 
The measure $\mu^d_R$ is obtained by keeping only the point masses 
of $\mu^d$ contained in the cube $\left[-\sqrt{d} , R + \sqrt{d}\right]^d$. This part of the proof (save for a small modification) already appears  in
\cite[Theorem 6]{MS}. 

Let $f$ be an integrable function and $\nu$ a finite sum of Dirac deltas.
By discretization, any lower bound for $c$ in 
$\alpha \left|\{ Mf \ge \alpha\}\right| \le c \|f\|_1$ is a lower
bound for $C$ in
$
\alpha \left|\{ M\nu \ge \alpha\}\right| \le C \nu (\mathbb{R}^d)
$, and viceversa. Here $\nu (\mathbb{R}^d)$ simply counts
the number of point masses in $\nu$. Observe that the cubes
used in Claim 4 to estimate the size of $M\mu^d(x)$, for $x\in[0,1]^d$,
never exceed a sidelength of $2\sqrt{d}$. 
Let
$\mu^d_R := \sum _{i} \delta _{x_i}$, where  
$R= R(d) >> d$ is a
  natural number
and $x_i\in \mathbb{Z}^d$ ranges over all the points with integer coordinates in the cube 
$[-\sqrt{d} , R + \sqrt{d}]^d$. Using the fact that the  estimates in the
preceding claims hold
for every unit subcube of $[0,R]^d$ with vertices in $\mathbb{Z}^d$
(by the same argument presented for $[0,1]^d$) 
we have
\begin{equation*}
c_d \ge \sup_{R > 0}  \frac{2^{-1}e^{t^2/2} \left|\left\{M_d\mu^d_R > 2^{-1}e^{t^2/2} \right\}\cap \left[0,R\right]^d\right|}{\mu^d_R (\mathbb{R}^d)} \end{equation*}
\begin{equation*}
 \ge \sup_{R > 0}  \left(\frac{R^d}{(R+2 \sqrt{d} + 1)^d}\right)
\left(\frac{t^{1/3} e^{-t^2/2} 2^{-1}e^{t^2/2}}{20 e^2 \sqrt{2\pi}} \right)
=
\frac{t^{1/3}}{40 e^2 \sqrt{2\pi}}.
\end{equation*}
\qed

\end{document}